\documentstyle[amssymb,amsfonts]{amsart}

\newtheorem{theorem}{Theorem}[section]
\newtheorem{corollary}[theorem]{Corollary}
\newtheorem{lemma}[theorem]{Lemma}

\newtheorem{conjecture}[theorem]{Conjecture}
\newtheorem{problem}[theorem]{Problem}
\newtheorem{question}[theorem]{Question}

\theoremstyle{definition}

\newtheorem{remark}[theorem]{Remark}

\def\R{{\hbox{\bf R}}}
\def\C{{\hbox{\bf C}}}

\def\P{{\hbox{\bf P}}}
\def\E{{\hbox{\bf E}}}

\font \roman = cmr10 at 10 true pt

\def\sign{{\hbox{ sign}}}
\def\Per{{\hbox{ Per}}}
\def\bv{{\hbox{\bf v}}}

\def\beq{\begin{equation}}
\def\eeq{\end{equation}}
\def\barray{\begin{eqnarray*}}
\def\earray{\end{eqnarray*}}

\def\bas{\begin{align*}}
\def\eas{\end{align*}}
\def\bi{\begin{itemize}}
\def\ei{\end{itemize}}

\def\corank{{\hbox{\roman corank}}}

\def\dist{{\hbox{\roman dist}}}

\def\Bp{{b^\prime}}

\def\emph#1{{\it #1}}
\def\textbf#1{{\bf #1}}

\def\Bp{{\mathbf p}}

\def\BZ{{\mathbf Z}}

\def\ep{{\epsilon}}

\theoremstyle{plain}

\title[Combinatorial problems in random matrix theory ]{Recent progress in combinatorial random matrix theory }

\author[Van H. Vu]
{Van H. Vu}

\begin{document}

\begin{abstract} 
	
	 We are going to discuss recent progress on many problems in random matrix theory of a combinatorial nature, including several breakthroughs that solve long standing famous 
	 conjectures. 
	
\end{abstract}
\maketitle

\section{Introduction}
The theory of random matrices is a  very rich topic in  mathematics. Besides being 
interesting in its  own right, random matrices play a 
fundamental role in various areas such as statistics, mathematical physics, combinatorics,  theoretical computer science,  etc. 
 
 In this survey, we focus on problems of a combinatorial nature. This is a continuation of an earlier suvery \cite{Vu14}. Since the publication of \cite{Vu14}, 
 the field has been flourishing and  there has been  significant  progress,  including the solutions of several major conjectures. New methods have been introduced which enable one to consider problems which seemed impenetrable only few years ago. Thus, an update is definitely in order.

  We mostly focus on the  following discrete models, noticing that the techniques developed for them usually work for a much wider class of ensembles. 
 
 \begin{itemize}
 
 \item $M_n$: random matrix of size $n$ whose entries 
 are i.i.d. Rademacher random variables (taking values $\pm 1$ with probability $1/2$). 
 In various papers, this is referred  as the random sign matrix or  Bernoulli matrices.

 \item $M_n^{sym}$: random 
 symmetric matrix of size $n$ whose (upper triangular) entries 
 are i.i.d. Rademacher  random variables. 
  
 \item Adjacency matrix of a random graph. 
 This matrix is symmetric and at position $ij$ we write $1$ if  $ij$  is an edge and zero otherwise.  
  
 \item Laplacian  of a random graph. 
 \end{itemize}

{\it Model of random graphs.}  We consider  two models: Erd\" os-R\'enyi and 
random regular graphs. For more information about these models, see 
\cite{Bolbook, AS, Wormald}.

\begin{itemize}

\item (Erd\"os-R\'enyi) We denote by $G(n,p)$ 
a random graph on $n$ vertices, generated by 
drawing an edge between any two vertices with probability $p$, independently. 

\item(Random regular graph) A random regular graph
on $n$ vertices with degree $d$ is obtained by sampling uniformly
over the set of all simple $d$-regular graphs on the vertex set
$\{1, \dots, n \}$. We denote this graph by $G_{n,d}$.

\end{itemize} 

It is important to notice that the edges of $G_{n,d}$ are not independent. Because of this,
this model is usually harder to study, compared to $G(n,p)$.


We denote
by $A(n,p)$ ($L(n,p)$) the adjacency (laplacian) matrix of the Erd\"os-R\'enyi random
graph $G(n,p)$ and by $A_{n,d}$ ($L_{n,d}$)  the adjacency (laplacian) matrix of 
$G_{n,d}$, respectively.

{\it Notation.} In the whole paper, we assume that $n$ is large. The asymptotic
notation such as $o, O, \Theta$ is used under the assumption that $n \rightarrow \infty$.
We write $A \ll B$ if $A= o(B)$.  $c$ denotes a universal
constant. All logarithms have natural base, if not specified
otherwise.

\section{The singular probability}\label{section:singular}
\label{rank}

The most famous combinatorial problem concerning random matrices is perhaps 
the "singularity" problem. Let $p_n$ be the
probability that $M_n$ is singular.   Trivially, 
$$p_n \ge 2^{-n}, $$ as the RHS is the probability that the first two rows are equal. 

By choosing  any two rows (columns) and also replacing equal by equal up to sign, one can have a slightly better lower bound 

\begin{equation} \label{precise} p_n \ge (4-o(1)) { n \choose 2} 2^{-n}  = (\frac{1}{2} +o(1))^n . \end{equation}  

The following conjecture is folklore.

\begin{conjecture} [Singularity, non-symmetric] \label{conjecture:singularity}  $p_n =(\frac{1}{2} +o(1))^n . $ \end{conjecture}

One  can formulate  even more precise conjectures, based on the following belief

{\bf Phenomenon I.}  {\it The dominating reason for singularity of a random matrix  is the dependency  between a few rows/columns}. 

For instance, \eqref{precise} suggests 

\begin{conjecture} \label{stronger} $p_n =(2+o(1)) n^2 2^{-n} . $ \end{conjecture}

Examining the dependence between 3,4, 5 etc rows (columns) would lead to stronger conjectures with smaller error terms; see \cite{Arat}.

It is already non-trivial to prove that  $p_n=o(1)$.  This was first done by Koml\'os  \cite{Kom}
 in 1967.
 Later, Koml\'os (see \cite{Bolbook}) 
 found a new proof  which gave quantitative bound $p_n=O(n^{-1/2})$. The first exponential bound is due  to Kahn,
Koml\'os,  and Szemr\'edi \cite{KKS}, who  proved  $p(n) \le .999^n$. Their arguments were simplified by Tao
and Vu in \cite{TVdet}, resulting in a slightly better bound
$O(.958^n)$.  Shortly afterwards, Tao and Vu \cite{TVsing} combined the approach from \cite{KKS} with
the idea of inverse theorems   (see \cite[Chapter 7]{TVbook} or \cite{NVsurvey}  for surveys)  to obtained a  more
significant improvement $p(n) \le (3/4 +o(1) ) ^n$. With an additional twist,  Bourgain, Vu,  and Wood  \cite{BVW} improved the bound further to
 $p(n) \le (\frac{1}{\sqrt 2}+o(1))^n $. In a different direction, Rudelson and Vershynin reproved KKS result in  a stronger form involving a lower bound for the least singular value 
 \cite{RV1}

 The method from \cite{TVsing, BVW} enables  one to deduce bounds on $p(n)$ directly  from 
 simple trigonometrical estimates. For instance, the $3/4$-bound comes from the fact that 
 $$|\cos x| \le \frac{3}{4} + \frac{1}{4} \cos 2x, $$ while the 
 $1/\sqrt 2$ bound come from 
 $$|\cos x|^2 =\frac{1}{2} + \frac{1}{2} \cos 2x. $$

 In 2018, Tikhomirov proved  Conjecture  \ref{conjecture:singularity} \cite{Tikomirov18}.
 
 \begin{theorem} 
 	$p_n =(\frac{1}{2} +o(1))^n . $
 	\end{theorem}

Each of the above mentioned papers contain new, highly non-trivial, ideas, but the core  of the matter  is a phenomenon called anti-concentration. 
 Tikhomirov's proof combines previous  ideas  with  a new,  powerful,  double counting argument (which he referred  to as {\it inversion of randomness} ). This, via a sophisticated discretization procedure,  reduces matters to studying
  the anti-concentration properties of random walks with random coefficients. This is perhaps the main difference from previous works which 
  considered  random walks with deterministic coefficients. The proof in \cite{Tikomirov18} also provided a bound on the  least singular value, extending 
  the result from \cite{RV1}.

To conclude this section, let us discuss   a classical   anti-concentration result.
Let $\bv =\{v_1, \dots, v_n\}$ be a set of $n$ non-zero real numbers  and
$\xi_1, \dots, \xi_n$ be i.i.d random Rademacher  variables. Define
$S:= \sum_{i=1}^n \xi_i v_i$, $p_{\bv} (a) = \P (S=a)$,  and
$p_{\bv} = \sup_{a \in \BZ} p_{\bv } (a)$.

The problem of estimating $p_{\bv}$ came  from a paper of Littlewood and Offord 
in the 1940s \cite{LO}, as a key technical ingredient in their study of real roots of random polynomials.  Erd\"os, 
improving a result of Littlewood and Offord, proved the following theorem, which we will refer
  to as the Erd\"os-Littlewood-Offord small ball inequality; see \cite{NVsurvey} for more details. 

\begin{theorem} \label{theorem:LO}   Let $v_1, \dots, v_n$ be non-zero
numbers and $\xi_i$ be i.i.d Rademacher  random variables. Then
$$ p_{\bv} \le \frac{{n \choose
{\lfloor n/2 \rfloor}}}{2^n} = O(n^{-1/2} ). $$\end{theorem}

Theorem \ref{theorem:LO} is a classical result in combinatorics and
have many non-trivial extensions 
with far reaching consequences 
 (see \cite{ Hal, SSz, NVsurvey}, \cite[Chapter 7]{TVbook}  and the references therein).

To give the reader a feeling about how  anti-concentration  estimates  can be useful in estimating $p_n$, let us sketch 
 the proof of  $p_n =o(1) $. We build $M_n$ by adding one random row at a time. 
  Assume that the first $n-1$ rows are 
independent and form a hyperplane with normal vector $\bv=(v_1, \dots, v_n)$. Conditioned on these rows, the probability that 
$M_n$ is singular  is 
$$\P ( X \cdot \bv =0) = \P (\xi_1 v_1 + \dots + \xi_ n v_n =0) , $$ where 
$X=(\xi_1, \dots, \xi_n)$ is the last row. If we can prove that many of the $v_i$ are non-zero, then Theorem \ref{theorem:LO} can be used to finish the job. 
In order to obtain strong quantitative  bounds on $p_n$,  the  following phenomenon proves useful 

{\bf Phenomenon II.} {\it If $\P ( X \cdot \bv =0)$ is relatively large, then the coefficients $v_1, \dots, v_n$ posses a  strong 
additive structure. }

For more discussion on this phenomenon and anti-concentration in general, we refer to the survey  \cite{NVsurvey}. 

\begin{remark} While polishing this survey, I learned of two new remarkable results (just added to the arxiv). First, Irmatov \cite{Irmatov} announced a proof of Conjecture \ref{stronger}. His approach seems quite different from all previous ones. Second, Litvak and Tikhomirov \cite{Litvak} announced a solution for a variant of  Conjecture \ref{stronger} in the spare case, for a very wide range of sparsity. 
	
\end{remark} 
\section{The singular probability: symmetric case}\label{section:singularsym}

As an analogue to the problem of the last section, it is natural to raise the question of estimating 
 $p_{n}^{sym}$, the probability that the symmetric matrix $M_n^{sym}$ singular. 
 
 This problem was mentioned to the author by Kalai and  Linial (personal conversations) 
 around 2004. To my surprise, at that point, even the analogue of Komlos' 1967 result was not known. 
 According to Kalai and Linial,  the following conjecture was circulated by  Weiss  in the 1980s, although it is quite possible that Koml\'os 
 had thought about it earlier.

 
 \begin{conjecture} \label{conj:weiss} 
  $p_n^{sym} =o(1) $. \end{conjecture}

The main difficulty concerning $M^{sym}_n$ is that its rows are no longer independent. This independence plays  a critical role in all results discussed in the previous section.

In \cite{CTV}, Costello, Tao, and Vu   found a way to circumvent the dependency. 
They  build  
the symmetric matrix $M_n^{sym}$ corner to corner. 
In step $k$, one considers the top left sub matrix of size $k$.  The strategy, following an idea of Koml\'os \cite{Kom},  is to show  that with high probability, 
the co-rank of this matrix, as $k$ increases, behaves like the end point of a  biased  random walk on non-negative integers  which has a strong tendency to go to the left whenever possible. This leads to a confirmation of Weiss' conjecture.

 \begin{theorem} \label{theorem:CTV} $p_n^{sym} = o(1) $.  \end{theorem}

The key technical tool in the proof of Theorem \ref{theorem:CTV} is 
the following (quadratic) variant of Theorem \ref{theorem:LO}.

Let us consider the last step in the process when  the $(n-1) \times (n-1)$ submatrix $M_{n-1}^{sym}$ has been built. To obtain $M_n^{sym}$, we add a random  row $X=(\xi_1, \dots, \xi_n)$ and its transpose. Conditioning on $M_{n-1}^{sym}$, the
determinant of the resulting $n \times n$ matrix is 
$$\sum_{1\le i, j \le n-1} a_{ij} \xi_i \xi_j  + \det M_{n-1} ,$$ where $a_{ij}$ 
(up to the signs) are 
the cofactors  of $M_{n-1}$. If $M_n^{sym}$ is singular, then its determinant is $0$, which implies 
$$Q:= \sum_{1\le i, j \le n-1} a_{ij} \xi_i \xi_j  = -\det M_{n-1} .$$ This gives ground for applying an anti-concentration result for quadratic forms.

Motivated by the non-symmetric case, it is natural to conjecture

\begin{conjecture} [Singularity, symmetric] \label{conjecture:ranksym}  $p^{sym} _n = (1/2+o(1))^n .$  \end{conjecture}

We leave it to interested readers to formulate more precise conjectures based on Phenomenon I. 
 The concrete bound from \cite{CTV} is $n^{-1/8}$, which can be easily improved to $n^{-1/4}$. 
 Costello \cite{Cos} improved the bound  to $n^{-1/2 +\epsilon}$ and Nguyen \cite{Nguyen1} pushed it further   to $n^{-\omega (1) }$. Next, Vershynin 
 proved a bound of the form  $\exp( -n^{c} )$, for some small constant $c >0$ \cite{Ver1}. 
 In   \cite{Ferber}, Ferber and Jain showed that one can take $c=1/4$.  This was improved further to $c =1/2$ by  
 Campos, Mattos, Morris, and  Morrison \cite{CMMM}.

\section{Ranks and co-ranks}

The singular probability is  the probability that the random matrix has co-rank at least one. {\it What about  larger co-ranks ? } 
Let us use $p_{n,k}$ to denote the probability that $M_n$ has co-rank at least $k$. It is easy to show that 
\begin{equation} p_{n,k} \ge (\frac{1}{2} +o(1)) ^{kn}. \end{equation} 

It is tempting to conjecture that this bound is sharp for constants $k$. In \cite{KKS}, Kahn,
 Koml\'os, and Szemer\'edi showed

\begin{theorem} \label{theorem:KKSrank}
There is a function $\epsilon(k)$ tending to zero  with $k$ tending to infinity such that 
$$p_{n,k} \le  \epsilon (k)  ^{n} . $$
\end{theorem}

In  \cite{BVW}, Bourgain, Vu,  and Wood  consider a variant of $M_n$ where the first $l$ rows are fixed
and the next $n-l$ are random. Let $L$ be the submatrix defined by the first $l$ rows and denote the model by $M_{n}(L)$. 
It is clear that  $\corank M_n(L) \ge \corank L $. The authors of \cite{BVW} showed
(\cite[Theorem 1.4]{BVW})

\begin{theorem} \label{theorem:BVWrank}
There is a positive constant $c$ such that 
$$ \P( \corank M_n(L) > \corank L) \le (1-c)^n. $$
\end{theorem}

Let us go back to  the symmetric model $M_n^{sym}$ and view it from this  new angle,
exploiting a connection to  Erd\" os-R\'enyi random graph $G(n,1/2)$. One can see that 
$$M_n^{sym}  = 2A(n,1/2) - J_n , $$ 
where $J_n$ is the all-one matrix of size $n$. (Here 
 we allow $G(n,1/2)$ to have loops, so the diagonal entries of $A(n,1/2)$ can be one. If we fix all diagonal entries to be zero, the analysis does not change significantly.) 
 Since $J_n$ has rank one,  it  follows from Theorem \ref{theorem:CTV}  that 
 with probablity $1-o(1)$, $A(n,1/2)$ has corank at most one.

 One can reduce the co-rank to zero by a slightly trickier argument. Consider $M_{n+1}^{sym}$ instead of $M_n^{sym}$ and normalize so that its first row and column are all- negative one. Adding this matrix with $J_{n+1}$, we obtain a matrix of the form\\ 
\begin{equation*}
\begin{pmatrix}
0&0\\
0&M_n^{sym}+J_n
\end{pmatrix}
\end{equation*}

 Thus we conclude
 \begin{corollary}
 With probability $1-o(1)$, $\corank A(n,1/2)=0$. 
 \end{corollary}

From the random graph point of view, it is  natural to ask if this statement  holds
for a  different density $p$ and if there is a {\it threshold} phenomenon; see \cite{AS} or \cite{Bolbook} for the definition of threshold. 
It is clear that the adjacency matrix is singular if  the density $p$ is very small. 
Indeed,  if $p < (1-\ep)\log n/n$, for any constant $\epsilon >0$, then by the coupon collector theorem, $G(n,p)$ has, with high probability,  isolated vertices,  which
correspond to zero rows in the adjacency matrix.  Costello and Vu    \cite{CV} proved that $\log n/n$ is the right
threshold.

\begin{theorem} \label{theorem:CVrank2} For any constant $\ep>0$,
with probability $1-o(1)$, 
$$ \corank A(n,(1+\ep)\log n/n)=0. $$ \end{theorem}

Strengthening Theorem \ref{theorem:CVrank2},  Basak and Rudelson \cite{Basak2} showed that one can replace $(1+ \epsilon)  \log n/n $ by $\log n/n + \gamma (n)/n$ where $\gamma (n)$ is any function tending to infinity.  In this direction, the most satisfying result is by  Addario-Berry and Eslava  \cite{Ahitting}, who   proved the following hitting time version. We generate the random graph by adding random edges one by one (the next  random edge is uniformly chosen from the set of all available edges). Let $T$ be the first time when the graph has  no isolated vertices.

\begin{theorem} \label{theorem:CVrank3} 
	With probability $1-o(1)$, the graph is full rank at time $T$.  \end{theorem}

 For $p < \log n/n$, the co-rank of $A(n,p)$ is no longer zero and its behavior, as a   random variable,  is not entirely understood. For the case  when 
 $p= c \log n/n$ for some constant $ 0 < c<1$, Costello and Vu \cite{CV2} showed that with probability $1-o(1)$, the co-rank is determined by small subgraphs whose structure forces the rank to drop,  which is consistent with {\bf Phenomenon I}. For example,  
 
\begin{theorem}  For any constant $\epsilon >0$
and $ (1/2 +\epsilon) \log n/n < p < (1-\epsilon) \log n /n$, 
with probability $1-o(1)$, 
$\corank A(n, p )$ equals  the number of isolated vertices.   \end{theorem}

For a smaller  $p$, one needs to take into account other small structures such as {\it cherries} (a cherry is  a pair of vertices of degree one with a common neighbor; 
in the matrix, this subgraph forces  two identical rows).  The main result of \cite{CV2} gives a precise formula for the co-rank in term of these parameters. 

When  $p =c/n, c>1$, $G(n,p)$ consists 
of a giant component and many small components. It makes sense to focus on the giant component  which we denote by  $Giant (n,p)$.  Since $Giant (n,p)$ has cherries , the adjacency matrix of $Giant (n,p)$ is singular (with high probability). 
However, if we look at the $k$-core of $Giant (n,p)$, for a sufficiently  large $k$, it seems plausible that this subgraph has full rank. (A $k$-core of a graph $G$ is a maximal connected subgraph of $G$ in which all vertices have degree at least $k$.)

\begin{conjecture} [$k$-core] Let $c >1$ be a constant and set  $p= c/n$. There is a constant $k_0$ such that for all $ k \ge k_0$ the following holds. 
With probability $1-o(1)$, the adjacency matrix of the $k$-core of $Giant(n,p)$ is non-singular. 
\end{conjecture}

  Bordenave, Lelarge,  and  Salez \cite{Bor} proved the following asymptotic   result 

\begin{theorem} Consider $G(n, c/n)$ for some constant $c >0$.  Then with probability $(1-o(1)) n$, $$rank (A (n, c/n)) = (2 - q - e^{-cq} -  cq e^{-cq } +o(1)) n , $$  where 
$0 < q<1$ is the smallest solution of  $q = \exp (-c  \exp {-cq} ) $.

\end{theorem}

In \cite{Coja19}, Coja-Oghlan,  Erg\"ur, Gao, Hetterich, and Rolvien studied random matrices with prescribed number of non-zeroes in each row and column and achieved an asymptotically sharp estimate for the rank; see  \cite{Coja19} for details.

\section {Random regular graphs}


Let  us  consider the random regular graph $G_{n,d}$. For $d=2$, $G_{n,d}$ is just the union of disjoint circles. It is not hard to show that with probability $1-o(1)$, one of these circles has length divisible by 4, and thus its adjacency matrix is singular (in fact, its corank is 
$\Theta (n)$ as the number of circles of length divisible by 4 is of this order). 
In \cite{Vu1}, the author raised the following conjecture  (which later appeared in \cite{Frieze14, Vu14} as well)

\begin{conjecture} [Singularity of Random regular graphs]\label{regular} For any $3 \le d \le n-1$, with probability $1-o(1)$
$A_{n,d}$ is non-singular. \end{conjecture}

Many of the earlier works on this conjecture  considered  the non-symmetric model, namely, random directed regular graphs. 
In this model, the matrix is chosen uniformly among all (not necessarily symmetric)$(0,1)$ matrices with  exactly $d$ ones in each column and row. For this non-symmetric model, 
Conjecture  \ref{regular} was confirmed for $C \ln^2 n \le d \le n - C \ln ^2 n$ by Cook \cite{Cook1}, and for $C_n \le d  \le n / C \ln^2 n  $ by  Litvak,  Lytova,  Tikhomirov,  Tomczak-Jaegermann, and  Youssef
\cite{Litvakregular}, where $C$ is a sufficiently large constant and $C_n$ tends to infinity arbitrarily slowly; see also \cite{Litvakregular3}  for an estimate on the least singular value. Furthermore, Litvak et al. \cite{Litvakregular2} also showed (whp)  that the rank is at least $n-1$ for any $d$. Huang showed that the rank equals $n$ (whp) for any fixed $d$; see \cite{Huang2}. 

For the (original) symmetric case. Landon, Sose, and Yau \cite{Landon}  showed that the conjecture holds for $d \ge n^c$ for any constant $c$, as a corollary of a more general
and precise  universality theorem. The most challenging case, $d$ being a constant, was solved recently by Meszaros \cite{Meszaros}  and Huang \cite{Huang1}. In particular, Huang proved

\begin{theorem}
	For any fixed $d \ge 3$, the probability that $A_{n,d}$ is singular is at most $n^{-c}$ for some constant $c >0$.
	\end{theorem}

	Huang's proof showed that one can take $c= \min \{1/8, (d-2)/(5d-6) \}$. On the other hand, he noted that the probability that $A_{n,d} $ is singular is at least 
	$n^{-d+2} $. It is interesting open question to find the sharp value of the exponent.

\section {Determinant and Permanent} \label{section:determinant} 

Let us  start with a basic question

\begin{question}  How big is the determinant of $M_n$?  \end{question} 

This was the original motivation of  Koml\'os' study, which started the line of research  discussed in  Section \ref{section:singular}; see  \cite{Kom, Kom2}.
However, the fact that $M_n$ is non-singular (whp) alone 
does not give any non-trivial estimate on the order of magnitude of $|\det M_n|$.

As  all rows of $M_n$ has length $\sqrt n$,  Hadamard's inequality 
implies that $| \det M_n | \le
n^{n/2}$. It  was   conjectured that with probability close to 1,
$|\det M_n|$ is close to this upper bound.

\begin{conjecture} \label{conj:determinant} Almost surely $| \det M_n | = n^{(1/2-o(1)) n}$.   \end{conjecture}

This conjecture is supported by a well-known observation of
Tur\'an.

\begin{equation}  \label{fact:Turan}     \E ((\det M_n)^2) = n!.   \end{equation} 
To verify this, notice that 
$$(\det M_n)^2 = \sum_{\pi, \sigma \in S_n} (-1) ^{\sign \pi +\sign \sigma}
\prod_{i=1}^n \xi_{i \pi(i) } \xi_{i \sigma(i)} . $$
By linearity of expectation  and the fact that $\E (\xi_i) =0$, we have 
$$\E (\det M_n)^2 =\sum_{\pi \in S_n} 1 = n! . $$
It follows immediately by Markov's bound that for any function
$\omega (n)$ tending to infinity with $n$,
$$|\det M_n | \le \omega (n)\sqrt {n!} ,$$ with probability tending to 1. 

A  statement  of  Girko  (the main result of \cite{Girko1, Girko2})  implies that 
$|\det M_n |$ is typically close to $\sqrt {n!} $. However, his proof appears to contain some 
gaps (see \cite{NVdet} for details).

\noindent In \cite{TVdet},  Tao  and Vu  established the
matching lower bound, confirming Conjecture
\ref{conj:determinant}.

\begin{theorem} \label{theorem:TVdet} With probability $1-o(1)$, 
$$|\det M_n |  \ge \sqrt {n!} \exp (- 29 \sqrt {n \log n}). $$ \end{theorem}

We sketch the proof very briefly as it contains a
useful lemma. 

First view  $|\det M_n|$ as the volume of the
parallelepiped spanned by $n$ random $\{-1,1\}$ vectors. This
volume is the product of the distances from the $(d+1)$st vector
to the subspace spanned by the first $d$ vectors, where $d$ runs
from $0$ to $n-1$. We are able to obtain a very tight control on
this distance (as a random variable), thanks to the following
lemma, which can be proved using  a powerful concentration  inequality by Talagrand  \cite{TVdet}.

\begin{lemma} \label{lemma:distance}
Let $W$ be a fixed subspace of dimension $1 \le d \le n-4$ and $X$
a random $\pm 1$ vector. For any $t >0$
\begin{equation}\label{pdw}
\P( |\dist(X, W) - \sqrt{n-d}| \ge t+1 ) \le  4 \exp(  -t^2/16).
\end{equation}
\end{lemma}

  The lemma, however, is not applicable
when $d$ is very close to $n$.  In this case, we need to make use of the fact  that $W$ is random.
Lemma \ref{lemma:distance}  appears handy  in many studies involving high dimensional probability.

Now we turn to the symmetric model $M_n^{sym}$. Again, by Hadamard's inequality 
$|\det M_n^{sym}| \le n^{n/2}$. 

\begin{conjecture} \label{conj:determinantsym} 
With probability $1-o(1)$  $$| \det M_n^{sym} | = n^{(1/2-o(1)) n}. $$  \end{conjecture}

Tur\'an's identity no longer holds
because of a correlation caused by symmetry. However, one can still show 
$$\E (\det M_n^{sym}) ^2 = n^{ (1+o(1)) n} . $$

On the other hand, proving a lower bound for $|\det M_n^{sym} |$
was  more difficult.  Recall that one  can interpret the determinant as the product of singular values, which, in this case, are the absolute values of the eigenvalues. By Wigner semi-circle  law (see \cite{Baibook}), we know the asymptotic of 
most of the singular values.  However, this law does not say anything about the smallest singular value which, in principle.  could be very close to zero. The problem of bounding 
the least singular value from below  was solved by  by  Nguyen \cite{Nguyen2} and Vershynin \cite{Ver1}.
 The results by Nguyen and Vershynin, combining with Wigner semi-circle law,  confirm Conjecture \ref{conj:determinantsym}

\begin{theorem} \label{theorem:determinantsym} 
With probability $1-o(1)$  $$| \det M_n^{sym} | = n^{(1/2-o(1)) n}. $$  \end{theorem}

The interested readers can find more precise results, which write down the limiting law of the determinant in 
\cite{Girko1, Girko2, NVdet, TVdetlaw, BMdet}. Let us now go  back to the random regular graphs. Consider the most interesting case when $d$ is a constant, we know that the matrix $A_{n,d} $ (whp) has full rank, so the determinant is non-zero. Its magnitude, however, is unknown. 

\begin{question} \label{conj:determinantsym} 
Estimate the determinant of $A_{n,d}$. Find the limiting law.   \end{question}

A completely open problem is to bound  the probability that $\det M_n$ (or other random determinant) takes on  a particular non-zero value. In contrast to the singularity conjecture 
(which addresses the case $\det M_n =0$), it seems that for any value $x \neq 0$, $\P (\det M_n =x)$ is sub-exponential. (This was first suggested to the author by Kalai in the early 2000s.) In fact, in view of Turan's identity, we conjecture that 

\begin{conjecture}  [Determinant] \label{det2} 
	For any $x \neq 0$, $\P (\det M_n =x) \le n^{-(1/2+o(1)) n }$. 
	
\end{conjecture}

The best current upper bound is exponential \cite{KKS}. A much weaker conjecture is that size of  the support of $\det M_n$ is super-exponential. 
But even this is not known.  The new developments  discussed in Section \ref{section:sandpile} may shed some light on this problem.

Let us now turn to the  permanent. 
Recall the  formal definition of the determinant of a matrix $M$ (with entries $m_{ij} , 1\le i, j \le n$)  
$$\det M: =\sum_{\pi \in S_n} (-1) ^{\sign \pi} \prod_{i=1}^n m_{i \pi (i) } . $$
The permanent of $M$ is defined as 
\begin{equation} \label{def:per} \Per M := \sum_{\pi \in S_n} \prod_{i=1}^n m_{i \pi (i) } . \end{equation}

Any question for determinant has its natural analogue for permanent. But typically, the problem becomes much harder as  permanent, unlike determinant,  does not admit  any useful 
geometric or linear algebraic interpretation. On the other  hand, it is easy to see that Tur\'an's identity still holds, namely
$$\E (\Per M_n)^2  = n ! . $$

\noindent It suggests that  $|\Per M_n|$ is typically $n^{(1/2-o(1))n }$. 
However, this was  much harder to prove than its determinant counterpart. The following conjecture, which is the  the permanent analogue  of Koml\'os' 1967  result $p_n=o(1)$,  was open for several decades

\begin{conjecture}
$\P (\Per M_n =0) =o(1)$. 
\end{conjecture}

In 2007, Tao and Vu \cite{TVper}  found  an entirely combinatorial approach to treat  the permanent problem, relying  on 
the formal definition \eqref{def:per} and making heavy use of martingale techniques from probabilistic combinatorics.   
They proved 

\begin{theorem} \label{theorem:permanent} 
With probability $1-o(1)$  $$|\Per M_n | = n^{(1/2-o(1)) n}. $$  \end{theorem}

As far as order of magnitude is concerned, the  still  missing piece of the picture is  the symmetric counterpart of 
Theorem \ref{theorem:permanent}. 

\begin{conjecture} \label{conj: permanentsym} 
With probability $1-o(1)$  $$|\Per M_n^{sym} | = n^{(1/2-o(1)) n}. $$  \end{conjecture}

Motivated by the singularity problem, it is 
 of interest  to find a strong estimate for the probability that the permanent is zero. In this aspect, we believe that determinant and permanent behave {\it differently}  and conjecture 
 
 \begin{conjecture} [Permanent]
 	The probability that $\Per M_n = 0$ is super exponentially small in $n$. 
 	\end{conjecture}
 	
 	The current bound is polynomial in $n$ \cite{TVper}. 
We do not know anything about the distribution of $|\Per M_n|$, either.  Even  simulation is challenging, as computing 
 permanent is a well known $\# P$ -complete problem; see  \cite{Valiant79}.

\section{ Graph expansion and the second eigenvalue}

Let $G$ be a  connected graph on $n$ points and $A$  its adjacency matrix with eigenvalues
 $\lambda_1 \ge \lambda_2 \ge \dots \ge \lambda_n$. If $G$ is $d$-regular  then $\lambda_1=d$ and by Perron-Frobenius theorem 
 no other eigenvalue has  larger absolute value. A parameter of  fundamental interest is 
 
 $$\lambda (G):= \max_{|\lambda_i | < d}   |\lambda_i | . $$

One can derive  interesting properties of the graph from the
value this parameter. The general phenomenon here is

 {\bf Phenomenon III.} {\it If  $\lambda(G)$ is significantly less than $d$, then the edges of $G$
distribute like in a  random graph with edge density $d/n$. }

A representative result  is the following
\cite{AS}. Let $A,B$ be sets of vertices and
$E(A,B)$ the number of edges with one end point in $A$ and the other in $B$, then 

\begin{equation} \label{lambda2} | E(A,B) -\frac{d}{n}|A||B| | \le \lambda(G) \sqrt{|A||B|} . \end{equation} 

Notice that the term $ \frac{d}{n}|A||B| $ is the expectation of the number of edges between $A$ and $B$ if
 $G$ is random (in the Erd\" os-R\'enyi sense) with edge density $d/n$. Graphs with small $\lambda$ are often called {\it pseudo-random} \cite{AS, CGW, KSsur}.

One can use this  information about edge
distribution  to derive various properties of the graph (see
\cite{KSsur} for many results of this kind). The whole concept can
be generalized for non-regular graphs, using  the Laplacian rather than the adjacency matrix (see, for
example, \cite{Chung}).

From \eqref{lambda2}, it is clear that the smaller $\lambda$, the more "random" is $G$.
{\it But how small can $\lambda$ be ? } 

In what follows, we restrict ourselves to  the most interesting case when $d$ is fixed and $n$ tends to infinity. 
In this case, Alon and Boppana (see \cite{Alon}) proved that   
$$\lambda(G) \ge 2 \sqrt{d-1} -o(1). $$
Graphs which satisfy $\lambda(G) < 2 \sqrt {d-1}$ are called Ramanujan graphs. 
It is very hard to construct such graphs explicitly, and all known constructions, such as those by 
Lubotzky, Phillip, and Sarnak \cite{LPS} and Margulis \cite{Mar} rely heavily on 
number theoretic results, which, in turn,  requires $d$ to have specific values.  A  more combinatorial  approach was found few years ago 
by Markus, Spielman, and Snivastava \cite{MSS}. Their  method (at least in the bipartite case) works for all $d$, but  the construction is not explicit.

\begin{theorem} 
\label{MSS} A bipartite Ramanujan graph exists for all fixed  degrees $d \ge 3$ and  sufficiently large $n$. 
\end{theorem}

While showing the existence of Ramanujan graphs is already highly non-trivial, the real question, in our opinion,  is to compute the limiting distribution of 
 $\lambda (G_{n,d})$,  which would lead to the exact probability of a random regular graph being 
 Ramanujan; see \cite{Miller} for a discussion and some numerical simulation.
 
 
A weaker conjecture, by 
Alon \cite{Alon},  asserts that for any fixed $d$, with probability $1-o(1)$ 
$$ \lambda_2 (G_{n,d}) = 2 \sqrt{d-1} +o(1). $$

\noindent Friedman \cite{Fri1} and Kahn and Szemer\'edi \cite{KSz} showed that if $d$ is fixed and $n$ tends to infinity, then with probability $1-o(1)$, 
$\lambda (G_{n,d}) = O(\sqrt d)$.  Friedman, in a highly
technical paper \cite{Fri2}, used the moment method to prove   Alon's  conjecture (see also \cite{FKoh} for a recent generalization)


\begin{theorem} \label{theorem:Friedman} 
For any fixed $d$ and $n$ tending to infinity,  with probability $1-o(1)$ 
$$\lambda (G_{n,d} ) =  2 \sqrt {d-1} +o(1).$$
\end{theorem}


For more recent developments concerning Friedman's theorem, including a new, shorter proof by Bordenave, see \cite{Puder, Bordenave}.


\section {Simple spectrum}

A matrix has simple spectrum if its eigenvalues are different. We discuss the following basic question

\begin{question}  Are random matrices simple ? \end{question} 

 It is easy to see that if the entries have  continuous distribution, then the spectrum is simple with probability 1. On the other hand, the discrete case is far from 
trivial.  In particular, Babai raised  the following conjecture in the 1980s:

\begin{conjecture}  \label{conjecture:Babai} 
	With probability $1-o(1)$, $G(n, 1/2)$ has a simple spectrum. 
\end{conjecture}

The motivation came from the well-known result (proved by Leighton-Miller and   Babai-Grigoriev-Mount; see \cite{Babai}) that the notorious  graph isomorphism problem is in {\bf P}  within the class of graphs with simple spectrum.  Few years ago,   Tao and Vu \cite{TVsimple} proved  this conjecture.

\begin{theorem}[Simple Spectrum] \label{simple}   Babai's conjecture holds. 
\end{theorem}

The same proof applies  for $M_n^{sym}$ (and many other ensembles). Let $s_n$ be the probability that the spectrum of $M_n^{sym}  $ is not simple. 
We observe that $s_n \ge 4^{-n} $, which is  the probability that 
the first 3 rows are the same (which guarantees  that zero is an eigenvalue with multiplicity at least 2). We conjecture

\begin{conjecture} [Simplicity] \label{simple-conj} $s_n = (4+o(1))^{-n} $.
	
	\end{conjecture}

The current best upper bound is  $s_n \le e^{-n^c}$ for some small constant $c>0$ \cite{NTVsimple}. 
Let us now formulate the singular value version of Babai's conjecture.

\begin{conjecture}  \label{simple-sing} 
	With probability $1-o(1)$, the singular values of $M_n^{sym} $ are different. 
\end{conjecture} 

Notice that the singular values of a symmetric matrix are the absolute values of its eigenvalues. Thus, this conjecture asserts that there is no two eigenvalues 
adding up to zero. 

One can pose the same questions for $M_n$.  In this direction, Ge \cite{Ge} proved the analogue of Theorem \ref{simple}, showing that 
with probability $1 -o(1) $, the spectrum of $M_n$ is simple. In a very recent paper, Luh and O'rourke \cite{LuhS} proved the first exponential bound, showing that 
 the probability that the spectrum of  $M_n$ is not simple  is at most $2^{-cn } $,  for some constant $c >0$. It looks plausible that Conjecture \ref{simple-conj} holds for $M_n$ as well. The $M_n$  analogue of Conjecture \ref{simple-sing}  is   also open.

\section{Normality } 

Another  basic notion in linear algebra is  that of {\it normality}. An $n \times n$ real matrix $A$  {\it normal} if  $AA^T=A^TA$. Few years ago, the author raised the following question.

\begin{question}
	How often is a random matrix normal?
\end{question}

  Despite the central role of normal matrices in matrix theory, to our surprise, we found no previous results concerning this  question. 
We consider  $M_n$ and denote by $\nu_n$ the probability that $M_n$ is normal. 
Clearly, the probability that $M_n$ is symmetric is $2^{-\left(0.5+o(1)\right)n^{2}}$. Since symmetric matrices are normal, 
$$ \nu_n  \ge 2^{-\left(0.5+o(1)\right)n^{2}}.  $$

We conjecture that this lower bound is sharp.

\begin{conjecture} [Normality]

	$$ \nu_n =  2^{-\left(0.5+o(1)\right)n^{2} }. $$
\end{conjecture}


In \cite{DVnormal}, Deneanu and Vu  proved

\begin{theorem}\label{thm1}
	
	$$\nu_n \le 2^{-\left(0.302+o(1)\right)n^{2}}.$$
\end{theorem}

Notice that in previous sections, the conjectural bounds are often of the form
$2^{-(c+o(1)) n}$, for some constant  $c >0$. While this probability is small, it is still much larger than $2^{- \Omega (n^2)} $, which enables one to exclude very rare events (those occurring with probability $2^{- \omega (n) }$) and then condition on their complement.

The difficulty with the normality problem is that we are aiming at a  bound which is extremely small (notice that any non-trivial event concerning $M_n$ holds  with probability {\it at least } $2^{-n^2}$, which is the mass of a single $\pm 1$ matrix).  There is simply no non-trivial event of probability $1 -2^{-\omega (n^2)}$ to condition on. Key to  \cite{DVnormal}   is a new  observation that for any given matrix, we can permute its rows and columns so that the ranks of certain submatrices follow a given pattern. The fact that there are only  $n! = 2^{o(n^2)} $ permutations works in our favor and enables us  to execute a different type of conditioning.

Ferber, Jain,  and Zhao \cite{FJZ} noticed that one lemma in \cite{DVnormal} can be improved, and reworking the whole argument one could improve the constant $.302$ slightly (maybe at the 5th decimal place). 

Another problem where $2^{- (.5+o(1)) n^2} $  could  be the right answer is bounding 
the probability that all eigenvalues of $M_n^{sym} $  are integers (they are apparently real).

\begin{conjecture} [Integral spectrum]
	The probability that $M_n^{sym} $ has an integral spectrum is $2^{- (.5 +o(1))n^2} $. 
	\end{conjecture}

	In \cite{Alon2},  Ahmadi, Alon,  Blake,   and  Shparlinski showed that the probability that $A(n,1/2) $ has an integral spectrum is at most $2^{-n/400} $. Costello and Williams 
	\cite{KW}  improved this bound to $2^{-cn^{3/2} } $, for some constant  $c >0$. Their proof can be modified to yield the same result for $M_n^{sym}$.   We also conjecture the 
	$M_n$ analogue of the above conjecture to hold 
	
	\begin{conjecture} [Gaussian integral spectrum]
		The probability that  all eigenvalues of $M_n$ are gaussian integers  is $2^{- (1 +o(1))  n^2} $. 
\end{conjecture} 

\section{Sandpile groups of random graphs} \label{section:sandpile}

 Given a graph $G$, its Laplacian is defined as 
 
 $$L(G)= A(G)- D(G) $$ where $A(G)$ is the adjacency matrix and $D(G)$ is a diagonal matrix whose $i$th entry is the degree of the $i$th vertex. If $G$ is $d$-regular, then 
 $L(G)= A(G)- dI$. 
 
 Let $Z$ be the set of integer vectors in $\R^n$ whose coordinates sum up to zero. It is clear that the row vectors of $L(G)$ is a subset of $Z$. 
 Let $R$ be the abelian group consisting of integer linear combinations of these vectors. The group $S:= Z/R$ is called the sand pile group of $G$. It is easy to see that

 $$| S| = |\det L(G) |, $$ which equals the number of spanning trees of $G$, by Kirkhoff's theorem. 
 
 In \cite{Wood}, Wood studied the structure of $S$, when $G=G(n, p)$ for any fixed $p$ \footnote{The Woods in this section and Section \ref{section:singular} are M. M. Wood and P. M. Wood, respectively.}. First, it is shown that for any fixed finite abelian group  $H$
 
 $$\P( S= H ) = o(1) . $$
 
 A finer question is the following. Any finite abelian group $H$ is the direct product of its Sylow subgroups. Now fix a prime $\Bp $ and a $\Bp$-group $H$. {\it What is the chance that the $\Bp$-Sylow subgroup of $S$ equals  $H$ ? } Wood \cite{Wood} proved that this probability is asymptotically

 $$ \frac{f(H) }{ |H| | Aut(H)  | } \prod _{j=0}^{\infty}   (1- \Bp^{-2j -1}  ),$$ where 
 $f(H)$ is the number of bilinear, symmetric, perfect maps from $H \times H$ to $\C^{\ast} $. 
 (For a concrete formula for $f(H)$, see \cite{Wood}.)  Wood noted that this  is similar (but not quite the same) to a formula suggested by Cohen-Lenstra heuristics.

 If $|S|= | \det L(G) | $ is divisible by $\Bp$, then its $\Bp$-Sylow subgroup must be non-trivial. Thus, one can use the result to compute the probability that 
 $|S|$ is divisible by $\Bp$. 
 
 \begin{corollary}
 	Let $\Bp (n,p)$ be the probability that $\det L(n,p))$ is divisible by $\Bp$, then 
 	
 	$$\Bp(n,p)= (1+o(1)) (1 -\prod_{j\ge 0} (1- \Bp^{-2j-1 })). $$
 	
 \end{corollary}
 
 For instance, the probability that the number of spanning trees of $G(n,p)$ is even is $\approx .5806...$. Interestingly, this probability does not depend on the density $p$, as long as it is fixed. 
In \cite{Meszaros}, Meszaros extended Wood's theorem to random regular graphs, and used this result to prove the non-singularity of random regular graphs with fixed degrees. 
See also \cite{Wood2, CLKW, Mapple1, Mapple2} for  related results in this direction.

Let us go back to $M_n$, which defines a map from $\BZ^n$ to itself. As shown in Section \ref{section:singular}, this map is (whp) injective. 
But how often is it {\it surjective } ?  Notice that $M_n$ is surjective iff $|\det M_n| =1$, thus the probability of being surjective tends to zero with $n$ as discussed in  Section \ref{section:determinant}.

From this point of view, a recent result of Nguyen and Wood \cite{NW}  is quite surprising. Consider a 
$n \times (n+1)$ random matrix with iid Rademacher  entries. This matrix defines a map from $\BZ^{n+1} $ to $\BZ^n$. {\it What is the probability that this map is surjective ?}   Nguyen and Wood  showed that this probability is  $$ (1+o(1)) \prod_{k\ge 2}  \zeta(k)^{-1} \approx .4358, $$  which is {\it bounded away}  from both 0 and 1. (Here $\zeta$ is the zeta function.)  An important step in  the proof shows that   with respect to different primes $\Bp_1, \dots, \Bp_k$ , the distributions of a random  determinant 
mod $\Bp_i$ are approximately independent.

\section{Miscellany}

An interesting (and seemingly hard) conjecture is the following, which came up in the conversation between the author and  P. M.  Wood 
in 2009.  Later,  Babai informed us that  he made the same conjecture (unpublished) in the 1970s. 

\begin{conjecture} [Irreducibility]
With probability $1-o(1)$, the characteristic polynomial of $M_n$ is irreducible. 
\end{conjecture}

In \cite{Vu14}, the author raised the following conjecture, which asserts that spectra can be used as finger prints.

\begin{conjecture} [Finger Print] 
A $\pm 1$ matrix is determined by its spectrum if  no other $\pm 1$ matrices (not counting trivial permutations) have the same spectrum. 
Then  almost all $\pm 1$ matrices are determined by their spectrum. 
\end{conjecture} 

One can  raise the same question for $M_n^{sym}$ or $G(n,1/2)$. It is known that there are non-isomorphic co-spectral graphs. However, these should form a negligible part of the set of 
all graphs. 

The matrix $M_n$ is (whp) non-symmetric. Thus, there is no obvious reason for its  to have many real eigenvalues. (The oddity of $n$ would guarantee  one real eigenvalue, but nothing more.) 
The following conjecture is motivated by our  joint work with Tao in  \cite{TVuniv-iid}.

\begin{conjecture} [Real Eigenvalues]
$M_n$  has, with high probability, $\Theta (\sqrt n)$ real eigenvalues.
\end{conjecture}

Edelman, Kostlan,  and Shub \cite{EKS}  obtained a  formula for the expectation of the number of real eigenvalues for a 
gaussian matrix  (which is or order $\Theta (\sqrt n)$). In \cite{TVuniv-iid}, Tao and Vu proved that the same formula holds 
(in the asymptotic sense) for certain random matrices with entries $(0, \pm 1)$. However, we do not know anything about  $M_n$. As a matter of fact, even the following 
"first step" seems hard

\begin{problem} [Two real roots]
Prove that  $M_n$ has, with high probability, at least 2 real eigenvalues.
\end{problem}

The next problem bears some resemblance to the famous "rigidity" problem in 
computer science. 
Given a $\pm 1$ square  matrix $M$, we denote by $Res (M)$ the minimum
number of entries we need to switch (from $1$ to $-1$ and vice
versa) in order to make $M$ singular. Thus, $Res$ can be seen as the {\it resilience} of the matrix against an effort to reduce  its rank. 
It is easy to show that $Res (M_n)$ is, with high probability, at most $(1/2+o(1))n$.

\begin{conjecture} [Rank Resilience] \label{conjecture:effect1} With probability $1-o(1)$, $Res (M_n)= (1/2+o(1))n. $  \end{conjecture}

\noindent For a recent partial result, see \cite{Luhrigid}.  
A closely related question (motivated by the notion of
local resilience from \cite{SudV}) is the following. Call a
$\{-1,1\}$ ($n$ by $n$) matrix $M$ {\it stubborn} if all matrices
obtained by switching (from $1$ to $-1$ and vice versa) the
diagonal entries of $M$ are non-singular (there are $2^n$ such
matrices).

\begin{conjecture} [Local resilience] \label{conjecture:effect2}  With probability $1-o(1)$, $M_n$ is stubborn.   \end{conjecture}

{\it Acknowledgement.} The author would like to thank  NSF and AFORS  for their generous support and K. Luh, C. Koenig, L. Addario-Berry, H. Nguyen, J. Huang, A. Litvak, and K. Tikhomirov for useful comments.

\end{document}